\documentclass[12pt]{amsart}
\usepackage{amsmath,amsfonts,amstext,amsthm,epsfig}

\newtheorem{Theorem}{Theorem}
\newtheorem{Lemma}{Lemma}%[section]
\makeatletter
    \addtocounter{section}{0}
    
     \@addtoreset{equation}{section}
\makeatother
\def\R{\mathbb R}

\begin{document}

\title[Absence of rapidly decaying solutions]
{On the absence of rapidly decaying solutions for parabolic
operators whose coefficients are non-Lipschitz continuous in time}%
\author{Daniele Del Santo}\author{Martino Prizzi}

\address{Daniele Del Santo, Universit\`a di Trieste, Dipartimento di
Matematica e Informatica,
Via Valerio 12/1, 34127 Trieste, Italy}%
\email{delsanto@univ.trieste.it}%

\address{Martino Prizzi, Universit\`a di Trieste, Dipartimento di
Matematica e Informatica, Via Valerio 12/1, 34127 Trieste, Italy}%
\email{prizzi@dsm.univ.trieste.it}%
%\thanks{}

\date{\today}%
%\dedicatory{}%
%\commby{}%
% ----------------------------------------------------------------
\begin{abstract}
We find minimal regularity conditions on the
coefficients of a parabolic operator, ensuring that no nontrivial
solution tends to zero faster than any exponential.
\end{abstract}
\maketitle
% ----------------------------------------------------------------

\section {Introduction, statements and remarks}

Let $A$ be a nonnegative self-adjoint operator in a Hilbert space
$H$. Consider the Cauchy problem
\begin{equation}
\begin{cases}\frac{du}{dt}+Au=0\\u(0)=u_0\end{cases}\label{M1}
\end{equation}
The solution $u(t)$ can be represented in terms of the spectral
resolution $E_\lambda$ of $-A$ and it turns out that its
asymptotic behavior is like $e^{-\lambda_0 t}$, where $\lambda_0$
is the infimum of those values of $\lambda$ for which $E_\lambda
u_0=u_0$. It follows that no solution, except the trivial one, can
tend to zero faster than any exponential.
\par Peter Lax \cite{lax} considered nonautonomous
perturbations of (\ref{M1}) of the form
\begin{equation}
\begin{cases}\frac{du}{dt}+(A+K(t))u=0\\u(0)=u_0\end{cases}\label{M2}
\end{equation}
where $K(t)$ is a bounded linear operator. He proved that, if the
norm of $K(t)$ is sufficiently small, then again solutions of
(\ref{M2}), unless identically zero, do not tend to zero faster
than any exponential.
\par The question then arised naturally, whether a similar result
could hold even for perturbations which were not relatively
bounded with respect to $A$. In the years following, attention
focussed mainly on parabolic inequalities, written in integrated
form, like
\begin{equation}
\int|\partial_t u
-\sum_{ij}a_{ij}(t,x)\partial_{x_i}\partial_{x_j}u|^2\,dx \leq
C_1(t)\int| u |^2\,dx+C_2(t)\int \sum_{i}|\partial_{x_i}u|^2\,dx.
\label{M3}
\end{equation}
Several results (see e.g. \cite{CoLees, lees, oga, prot}) were
obtained, relating the decay of $C_1(t)$, $C_2(t)$ and $\|\nabla_x
a_{ij}(t,\cdot)\|_{L^\infty}$ to that of the solutions of
(\ref{M3}). Some years later,  Agmon and Nirenberg \cite{AN}
reconsidered the whole matter by an abstract point of view and
proved a general result for inequalities of the form
\begin{equation}
\|\frac{du}{dt}+A(t)u\|\leq \Phi(t)\|u\|\label{M4}
\end{equation}
in a Banach space $X$.
\par Without entering into technical details, we notice that there
is a common feature in all the above mentioned results: at a
certain point one needs to perform some integration by parts and
this requires some (kind of) differentiability of the coefficients
with respect to $t$. That a certain amount of regularity were
actually {\em necessary} in order to get lower bounds for the
solutions became clear thanks to a well known example of Miller
\cite{mill}. He exibited a parabolic operator whose coefficients
are H\"older continuous of order $1/6$ with respect to $t$ and
which possesses solutions vanishing within a finite time.
\par The aim of this paper is the following: for a parabolic
inequality of the form (\ref{M3}), find the minimal regularity of
the coefficients $a_{ij}$'s with respect to $t$, ensuring that no
solution, except the trivial one, can tend to zero faster than any
exponential. \par We prove that a sufficient regularity condition is
given in terms of a modulus of continuity
satisfying the so called {\em Osgood condition}. The counter
example contained in \cite{DSP} shows that this condition is
optimal. The main result (Theorem 1 below) is a consequence of a
{\em Carleman estimate} in which the weight function depends on
the modulus of continuity; such kind of weight functions in
Carleman estimates were introduced by Tarama \cite{tara} in the
study of second order elliptic operators.
\par In order to make the presentation simpler, we consider
an equation whose coefficients are independent of the space
variable $x$. The general case can be recovered by the same
microlocal approximation procedure exploited in \cite{DSP}.

\par Let $a$ be a continuous function defined on ${\R}^+$ such that
\begin{equation}
Ê\Lambda_0^{-1}\leq a(t)\leq \Lambda_0\label{1}
\end{equation}
for some $\Lambda_0\geq 1$ and for all $t\in {\R}^+$. Let $\varphi$ be a
positive
function in $L^1({\R}^+)$.
Let $u$ be a function defined on ${\R}^+_t\times \R_x$ such that
\begin{equation}
u\in L^2_{\rm loc}({\R}^+_t, H^2({\R}_x))\cap H^1_{\rm loc}({\R}^+_t,
L^2({\R}_x)) \label {2}
\end{equation}
and
\begin{equation}
\|u_t(t, \cdot)-a(t)u_{xx}(t, \cdot)\|^2_{L^2({\R}_x)}\leq
\varphi(t) \|u(t, \cdot)\|^2_{H^1({\R}_x)} \label {3}
\end{equation}
for a.e. $t\in {\R}^+$.
A function $u$ satisfying the conditions ({\ref{2}) and (\ref{3}) is called
{\it rapidly decaying solution} to
Ê(\ref{3}) if for all $\lambda >0$,
\begin{equation}
\lim_{t\to +\infty} e^{\lambda t}\|u(t, \cdot)\|_{H^1({\R}_x)}=0.\label {4}
\end{equation}

Let $\mu$ be a {\it modulus of continuity} i.e. $\mu$ is a function defined
on $\R^+$ with values in
$\R^+$ such that $\mu$ is continuous, increasing, concave and $\mu(0)=0$.
A modulus of continuity Ê$\mu$ is said to satisfy the {\it Osgood
condition} if
\begin{equation}
\int_0^1{1\over \mu(s)}\, ds=+\infty.\label {5}
\end{equation}

Now we can state our main result:
\begin{Theorem}
Let $\mu$ be a modulus of continuity satisfying the Osgood condition.
Suppose that there exists a positive function $\psi$ in $ L^1(\R^+)\cap
ÊL^\infty (\R^+)$ such that
\begin{equation}
\sup_{\max\{0, t-{1\over 2}\}< t_1< t_2< t+{1\over 2}}
{|a(t_2)-a(t_1)|\over \mu(t_2-t_1)}\leq \psi(t)\label{6}
\end{equation}
for a.e. $t\in {\R}^+$.

If $u$ is a rapidly decaying solution to (\ref{3}) then $u\equiv 0$.
\label{t1}
\end{Theorem}

The counter example alluded to above is given by the following
\begin{Theorem} Let $\mu$ be a modulus of continuity which does not satisfy the Osgood 
condition.
Then there exists $l\in C({\mathbb R}_t)$ with $1/2\leq l(t)\leq 3/2$
for all
$t\in {\mathbb R_t}$ and
\begin{equation}
\sup_{0< |t_1- t_2| <1\atop t_1,t_2\in\R_t}
{|l(t_2)-l(t_1)|\over \mu(t_2-t_1)}<\infty
\end{equation}
  and
there exists $u$, $b_1$, $b_2$, $c\in C^\infty_b({\mathbb R}_t
\times{\mathbb R}^2_x)$ with ${\rm supp}\;
u=\{t\leq 1\}$ such that
\begin{equation}
\partial_t u-( \partial^2_{x_1}u+ l
\partial^2_{x_2}u)+b_1\partial_{x_1}u+b_2\partial_{x_2}u+cu=
0\quad {\it in}\ {\mathbb R}_t\times{\mathbb R}^2_x.
\end{equation}
\end{Theorem}
The proof of Theorem 2 is contained in our previous paper \cite{DSP}.

\section {Proof of Theorem 1}

First of all we remark that it is not restrictive to suppose
that
$\int_{t_1}^{t_2} \varphi(s)\, ds>0$
and
$\int_{t_1}^{t_2} \psi(s)\, ds>0$ for all $0\leq t_1<t_2$.
Moreover we will admit without lack of generality that
\begin{equation}
\int_{0}^{1} \varphi(s)\, ds\geq 1. \label{2.0}
\end{equation}
Let $\alpha>0$. We set, for
$t\geq 0$,
\begin{equation}
b(t)=\exp(-\alpha \int_0^t \varphi (\eta)\, d\eta.)\label{2.1}
\end{equation}
Let $\nu$ be a function defined in $[1,+\infty[$ such that
\begin{equation}
\nu(t)=\int_{1/t}^1 {1\over \mu(s)}\, ds;\label{2.2}
\end{equation}
we remark that (\ref{5}) gives, in particular,
$\nu([1,+\infty[)=[0,+\infty[$.
For $\gamma>0$ and $\tau\geq 0$ we define
\begin{equation}
\Psi_\gamma(\tau)=\nu^{-1}(\gamma\int_0^{\tau/\gamma}\psi(s)\, ds).\label{2.3}
\end{equation}
Finally we set, for $\gamma>0$ and $t\geq 0$,
\begin{equation}
\Phi_\gamma(t)=\int_0^t\Psi_\gamma(\gamma\eta){1\over
b(\eta)}(\int_0^\eta b(s)\varphi(s)\, ds)\, d\eta.\label{2.4}
\end{equation}
\begin{Lemma}
For all $\alpha>0$ there exists $\gamma_0>0$ such that
\begin{equation}
\begin{array}{lll}
&\displaystyle{\int_1^{+\infty}
b(t)\,
e^{2\Phi_\gamma(t)}\|v_t(t,\cdot)-a(t)v_{xx}(t,\cdot)\|^2_{L^2(\R_x)}\,dt}\\[0.3
cm]
&\displaystyle{\quad\geq {\alpha\over \Lambda_0}\int_1^{+\infty}
b(t) \varphi(t)\,
e^{2\Phi_\gamma(t)}
\|v_x(t,\cdot)\|^2_{L^2(\R_x)}\, dt }\\[0.3 cm]
&\displaystyle{\qquad\quad+\int_1^{+\infty}\Psi_\gamma(\gamma
t)\,b(t) \varphi(t)\,
e^{2\Phi_\gamma(t)}
\|v(t,\cdot)\|^2_{L^2(\R_x)}\, dt}
\end{array}\label{2.5}
\end{equation}
for all $\gamma\geq \gamma_0$ and for all  $v\in L^2(]1,+\infty[, H^2({\R}_x))\cap
H^1(]1,+\infty[, L^2({\R}_x))$ with {\em compact support}. \label{l1}
\end{Lemma}

Let us show how to prove Theorem~\ref{t1}
from the Carleman estimate
(\ref{2.5}). Let $w$ be a function in $L^2_{\rm loc}(]1,+\infty[,
H^2({\R}_x))\cap
H^1_{\rm loc}(]1,+\infty[, L^2({\R}_x))$ such that $w(t,x)=0$ for all
$(t,x)\in [0,2]\times \R$.
Suppose that $w$ satisfies
\begin{equation}
\lim_{t\to +\infty} e^{\lambda t}\|w(t, \cdot)\|_{H^1({\R}_x)}=0. \label
{2.5bis}
\end{equation}
for all $\lambda>0$.
We show first that an inequality similar to (\ref{2.5}) holds for $w$.

Consider $\chi\in C^\infty(\R)$ with $\chi$ decreasing,
Ê$\chi(s)=1$ for $s\leq 1$ and $\chi(s)=0$ for $s\geq 2$ and
define $v_n(t,x)=\chi(t/n)w(t,x)$. Then $v_n\in L^2(]1,+\infty[,
H^2({\R}_x))\\\cap H^1(]1,+\infty[, L^2({\R}_x))$ and is compactly
supported, so that by (\ref{2.5}) we deduce $$
\begin{array}{lll}
&\displaystyle{2\int_1^{+\infty}
b(t)\, e^{2\Phi_\gamma(t)}\chi^2({t\over
n})\|w_t(t,\cdot)-a(t)w_{xx}(t,\cdot)\|_{L^2}^2\,dt}\\[0.3 cm]
&\displaystyle{\quad\geq
{\alpha\over \Lambda_0}\int_1^{+\infty} b(t) \varphi(t)\,
e^{2\Phi_\gamma(t)}\chi^2({t\over
n})\|w_x(t,\cdot)\|_{L^2}^2\, dt }\\[0.3 cm]
&\displaystyle{\quad\qquad+\int_1^{+\infty}\Psi_\gamma(\gamma
t)\,b(t) \varphi(t)\,
e^{2\Phi_\gamma(t)}\chi^2({t\over
n})\|w(t,\cdot)\|_{L^2}^2\, dt}\\[0.3 cm]
&\displaystyle{\qquad\qquad\quad-2\int_1^{+\infty} b(t)\,
e^{2\Phi_\gamma(t)}({1\over n}\chi'({t\over n}))^2\|w(t,\cdot)\|_{L^2}^2\,dt}
\end{array}
$$ for all $\gamma\geq \gamma_0$. Remark now that
$\Psi_\gamma(\gamma t)\leq \nu^{-1}(\gamma \|\psi\|_{L^1(\R_x)})$
for all $\gamma$ and $t$, while $e^{- \alpha
\|\varphi\|_{L^1(\R_x)}}\leq b(t)\leq 1$ for all $t$. Consequently
we have $$ \Phi_\gamma(t)\leq
Ê\nu^{-1}(\gamma\|\psi\|_{L^1})\|\varphi\|_{L^1} e^{\alpha
\|\psi\|_{L^1}}t = C_\gamma t $$ for all $\gamma$ and $t$. Hence,
using (\ref{2.5bis}) and the fact that $w\in C(]1,+\infty[,
H^1({\R}_x))$ (see \cite[pp. 18-19]{LM}), we deduce that $$ b(t)
\varphi(t) e^{2\Phi_\gamma(t)}\chi^2({t\over
n})\|w_x(t,\cdot)\|_{L^2}^2\leq K_\gamma \varphi(t), $$ $$
\Psi_\gamma(\gamma t)\,b(t) \varphi(t)\,
e^{2\Phi_\gamma(t)}\chi^2({t\over n})\|w(t,\cdot)\|_{L^2}^2\leq
K'_\gamma \varphi(t) $$ and $$ b(t)\, e^{2\Phi_\gamma(t)}({1\over
n}\chi'({t\over n}))^2\|w(t,\cdot)\|_{L^2}^2 \leq K''_\gamma
e^{-\tilde \lambda t} $$ for a.e. $t$. Passing to the limit for
$n\to +\infty$, and applying the dominated convergence theorem on
the right hand side and the monotone convergence theorem on the
left hand side, we obtain that
\begin{equation}
\begin{array}{lll}
&\displaystyle{\int_1^{+\infty}
b(t)\,
e^{2\Phi_\gamma(t)}\|w_t(t,\cdot)-a(t)w_{xx}(t,\cdot)\|_{L^2}^2\,dt}\\[0.3
cm] Ê&\displaystyle{\quad\geq {\alpha \over 2\Lambda_0}\int_1^{+\infty}
b(t) \varphi(t)\,
e^{2\Phi_\gamma(t)}\|w_x(t,\cdot)\|_{L^2}^2\, dt }\\[0.3 cm]
&\displaystyle{\quad\qquad+{1 \over
2}\int_1^{+\infty}\Psi_\gamma(\gamma t)\,b(t)
\varphi(t)\, e^{2\Phi_\gamma(t)}\|w(t,\cdot)\|_{L^2}^2\, dt}
\end{array}\label{2.6}
\end{equation}
for all $\gamma\geq \gamma_0$.

Let now $u$ be a rapidly decaying solution to (\ref{3}). Let $\theta\in
C^\infty(\R)$ with $\theta$ increasing, Ê$\theta(s)=0$ for $s\leq 2$ and
$\theta(s)=1$ for $s\geq 3$. Setting Ê$w(t,x)=\theta(t)u(t,x)$ and applying
(\ref{2.6}) we obtain
$$
\begin{array}{l}
\displaystyle{ \int_1^3 \! b(t)\, e^{2\Phi_\gamma(t)}\|(\theta
u)_t-a(t)(\theta u)_{xx}\|_{L^2}^2\, dt
%}\\[0.3 cm] \displaystyle{\quad
+\int_3^{+\infty} \! b(t)\,
e^{2\Phi_\gamma(t)}\|u_t-a(t)u_{xx}\|_{L^2}^2\, dt}\\[0.3
cm] \displaystyle{
%\qquad\quad
=\int_1^{+\infty} b(t)\,
e^{2\Phi_\gamma(t)}\|w_t-a(t)w_{xx}\|_{L^2}^2\,
dt}\\[0.3 cm] \displaystyle{
%\qquad\qquad\quad
\geq {\alpha\over
2\Lambda_0}\int_1^{+\infty} \!\!
b(t) \varphi(t)\, e^{2\Phi_\gamma(t)}\|w_x\|_{L^2}^2\, dt
%}\\[0.3 cm] \displaystyle{\qquad\qquad\qquad
+ {1\over 2}\int_1^{+\infty} \!\!
\Psi_\gamma(\gamma Êt)\, b(t)\varphi(t)\, e^{2\Phi_\gamma(t)}\|w\|_{L^2}^2\, dt
}\\[0.3 cm]
\displaystyle{
%\qquad\quad\qquad\qquad\quad
\geq {\alpha\over
2\Lambda_0}\int_3^{+\infty}\! b(t) \varphi(t)\,
e^{2\Phi_\gamma(t)}\|u_x\|_{L^2}^2\, dt
%}\\[0.3 cm] \displaystyle{\qquad\qquad\qquad\qquad\quad
+ {1\over
2}\int_3^{+\infty} \!
\Psi_\gamma(\gamma
t)\, b(t)\varphi(t)\, e^{2\Phi_\gamma(t)}\|u\|_{L^2}^2\, dt. }
\end{array}
$$
Hence, using also (\ref{3}) we have
$$
\begin{array}{l}
\displaystyle{ \int_1^3 b(t)\, e^{2\Phi_\gamma(t)}\|(\theta
u)_t-a(t)(\theta u)_{xx}\|_{L^2}^2\, dt}\\[0.3
cm]
\displaystyle{\quad\geq \int_3^{+\infty} b(t) ({\alpha\over
2\Lambda_0}-1)\varphi(t)\,
e^{2\Phi_\gamma(t)}\|u_x\|_{L^2}^2\, dt }\\[0.3 cm]
\displaystyle{\qquad\quad+ \int_3^{+\infty} ({1\over 2}\Psi_\gamma(\gamma
t)-1)\, b(t)\varphi(t)\, e^{2\Phi_\gamma(t)}\|u\|_{L^2}^2\, dt .}
\end{array}
$$ We take $\alpha=2\Lambda_0$. We recall that $b(t)\leq 1$ and
that $\Phi_\gamma$ is increasing. Hence $$ \int_1^3 b(t)\|(\theta
u)_t-a(t)(\theta u)_{xx}\|_{L^2}^2\, dt \geq
\int_3^{+\infty}({1\over 2}\Psi_\gamma(\gamma t )-1) \,
b(t)\varphi(t)\|u\|_{L^2}^2\, dt $$ for all $\gamma\geq \gamma_0$.
Since $\Psi_\gamma (\gamma t)\geq \Psi_\gamma (\gamma)$ for all
$t\geq 1$ we obtain $$ \int_3^{+\infty}({1\over 2}
\Psi_\gamma(\gamma t)-1)\, b(t)\varphi(t)\|u\|_{L^2}^2\, dt \geq
({1\over 2} \Psi_\gamma(\gamma)-1) \int_3^{+\infty}
b(t)\varphi(t)\|u\|_{L^2}^2\, dt. $$ From (\ref{5}) we deduce that
$\lim_{\gamma\to +\infty}\Psi_\gamma(\gamma)=+\infty$ and
consequently letting $\gamma$ go to $+\infty$ we obtain that
$u(x,t)=0$ in $[3, +\infty[\times \R$. We apply now the backward
uniqueness result in \cite{DSP} and we easily deduce that $u\equiv
0$.

Let us come to the proof of Lemma~\ref{l1}. Setting
$z(t,x)=e^{\Phi_\gamma(t)}v(t,x)$ we have
$$
\begin{array}{l}
\displaystyle{\int_1^{+\infty} b(t)\,e^{2\Phi_\gamma(t)}\|v_t(t,\cdot)
-a(t)v_{xx}(t,\cdot)\|^2_{L^2(\R_x)}\,dt}\\[0.3 cm]
\displaystyle{= \int_1^{+\infty} b(t)\|z_t(t,\cdot)
-a(t) z_{xx}(t,\cdot)-\Phi'_\gamma(t) z(t,\cdot)\|_{L^2}^2\, dt}\\[0.3 cm]
Ê\displaystyle{= \int_1^{+\infty}\!\int_{\R_\xi} b(t)|\hat z_t(t,\xi)|^2\,d\xi
dt +\! \int_1^{+\infty}\!\int_{\R_\xi} b(t)(a(t)\xi^2-\Phi'_\gamma(t))^2 |\hat
z(t,\xi)|^2\,d\xi dt}\\[0.3 cm]
\qquad\qquad\qquad\qquad\qquad \displaystyle{+2\Re\int_1^{+\infty}\!\int_{\R_\xi}
b(t)(a(t)\xi^2-\Phi'_\gamma(t)) Ê\hat z_t(t,\xi) \overline{\hat
z(t,\xi)}\,d\xi dt}
\end{array}
$$
where $\hat z$ denotes the Fourier transform of $z$ with respect to the $x$
variable.
We compute the second part of the last term of the above inequality and we
obtain
$$
\begin{array}{l}
\displaystyle{ -2\Re\int_1^{+\infty}\int_{\R_\xi}b(t)
\Phi'_\gamma(t)\hat z_t(t,\xi) \overline{\hat z(t,\xi)}\,d\xi
dt}\\[0.3 cm] \quad
\displaystyle{=\int_1^{+\infty}\Psi_\gamma(\gamma t) b(t)
\varphi(t)\|z(t,\cdot)\|^2_{L^2}\,dt }\\[0.3 cm]
\qquad\quad\displaystyle{+\int_1^{+\infty}\int_{\R_\xi}
\gamma\Psi_\gamma'(\gamma t)(\int_0^t b(s)\varphi(s)\,ds) |\hat
z(t,\xi)|^2\,d\xi dt.}
\end{array}
$$ It remains to estimate the quantity $$
2\Re\int_1^{+\infty}\int_{\R_\xi} b(t)a(t)\xi^2 Ê\hat z_t(t,\xi)
\overline{\hat z(t,\xi)}\,d\xi dt. $$ Since $a$ is not
Lipschitz-continuous and consequently we cannot integrate by
parts, we exploit the approximation technique developed in
\cite{DSP}. Let $\rho\in C_0^\infty(\R)$ with ${\rm supp}\,
\rho\subseteq [-1/2,1/2]$, $\int_{\R}\rho(s)\, ds=1$ and
$\rho(s)\geq 0$ for all $s\in\R$. We set $$
a_\varepsilon(t)=\int_{\R}a(s){1\over \varepsilon}\rho({t-s\over
\varepsilon})\, ds, $$ where $a$ has been extended to $\R$ setting
$a(t)=a(0)$ for all $t\leq 0$. We obtain that there exists $C_0>0$
such that $$ |a_\varepsilon(t)-a(t)|\leq \mu(\varepsilon)\,\psi(t)
$$ and $$ |a'_\varepsilon(t)|\leq C_0{\mu(\varepsilon)\over
\varepsilon}\, \psi(t) $$ for all $\varepsilon\in\, ]0,1]$ and for
a.e. $t\in \R^+$. Hence $$
\begin{array}{l}
\displaystyle{2\Re\int_1^{+\infty}\int_{\R_\xi}
b(t) a(t)\xi^2 \hat z_t(t,\xi)\overline{\hat
z(t,\xi)}\,d\xi dt}\\[0.3 cm]
\quad\displaystyle{=2\Re\int_1^{+\infty}\int_{\R_\xi}
b(t) a_\varepsilon(t)\xi^2 \hat z_t(t,\xi) \overline{\hat
z(t,\xi)}\,d\xi dt}\\[0.3 cm]
\qquad\quad\displaystyle{+2\Re\int_1^{+\infty}\int_{\R_\xi}
b(t) (a(t)-a_\varepsilon(t))\xi^2 \hat z_t(t,\xi) \overline{\hat
z(t,\xi)}\,d\xi dt.}
\end{array}
$$
We have
$$
\begin{array}{l}
\displaystyle{2\Re\int_1^{+\infty}\int_{\R_\xi}
b(t) a_\varepsilon(t)\xi^2 \hat z_t(t,\xi) \overline{\hat
z(t,\xi)}\,d\xi dt}\\[0.3 cm]
\displaystyle{\quad=-\int_1^{+\infty}\int_{\R_\xi}(b(t)a_\varepsilon(t))'
\xi^2 |\hat z(t,\xi)|^2\,d\xi dt}\\[0.3 cm]
\displaystyle{\qquad\geq \int_1^{+\infty}\int_{\R_\xi}
b(t)(\alpha \varphi(t) a_\varepsilon(t)-|a'_\varepsilon(t)|)
\xi^2 |\hat z(t,\xi)|^2\,d\xi dt}\\[0.3 cm]
\displaystyle{\qquad\quad\geq {\alpha\over \Lambda_0}\int_1^{+\infty}
b(t) \varphi(t) \|z_x(t,\cdot)\|_{L^2}^2\, dt}\\[0.3 cm]
\displaystyle{\qquad\qquad\qquad -C_0
\int_1^{+\infty}\int_{\R_\xi} b(t) \psi(t) {\mu(\varepsilon)\over \varepsilon}
\xi^2 |\hat z(t,\xi)|^2\,d\xi dt}
\end{array}
$$
and
$$
\begin{array}{l}
\displaystyle{2\Re\int_1^{+\infty}\int_{\R_\xi}
b(t) (a(t)-a_\varepsilon(t))\xi^2 Ê\hat z_t(t,\xi)\overline{\hat
z(t,\xi)}\,d\xi dt}\\[0.3 cm]
\displaystyle{\geq
-\int_1^{+\infty}\int_{\R_\xi} b(t)|\hat z_t(t,\xi)|^2\,d\xi dt
- \int_1^{+\infty}\int_{\R_\xi} Êb(t)\psi^2(t) \mu^2(\varepsilon)
\xi^4|\hat z (t,\xi)|^2\,d\xi dt.}
\end{array}
$$ Putting all these inequalities together it is easy to see that
(\ref{2.5}) will be a consequence of the following claim: \par for
all $\alpha>0$ there exist $\gamma_0>0$ and a function $\R\to
\,]0,1]$, $\xi\mapsto \varepsilon_\xi$ such that
\begin{equation}
\begin{array}{ll}
\displaystyle{\int_1^{+\infty}\int_{\R_\xi}
(b(t)(a(t)\xi^2-\Phi'_\gamma(t))^2 +\gamma\Psi_\gamma'(\gamma t)\int_0^t
b(s)\varphi(s)\,ds) |\hat z(t,\xi)|^2\, d\xi dt}\\[0.5 cm]
\qquad\displaystyle{-\int_1^{+\infty}\int_{\R_\xi} b(t) \psi(t)(C_0
{\mu(\varepsilon_\xi)\over \varepsilon_\xi}\xi^2+ \psi(t)
\mu^2(\varepsilon_\xi)\xi^4) |\hat z(t,\xi)|^2\,d\xi dt \geq 0}\label{2.7}
\end{array}
\end{equation}
for all $\gamma\geq\gamma_0$ and for all
$z(t,x)=e^{\Phi_\gamma(t)}v(t,x)$, provided $v\in L^2(]1,+\infty[,
H^2({\R}_x))\cap H^1(]1,+\infty[, L^2({\R}_x))$ is compactly supported. \par From
(\ref{2.2}) and (\ref{2.3}) we have that
\begin{equation}
\Psi'_\gamma(\gamma t) = \Psi^2_\gamma(\gamma t)\mu({1\over \Psi_\gamma(\gamma
t)})\psi(t).\label{2.9}
\end{equation}
The concavity of $\mu$ implies that the function $\sigma\mapsto\sigma
\mu(1/\sigma)$ is increasing on $[1,+\infty[$ and consequently the function
$\sigma\mapsto\sigma^2\mu(1/\sigma)$ is increasing and Ê$\sigma^2\mu(1/\sigma)
\geq \sigma \mu(1)$ for all $\sigma\in [1,+\infty[$. Hence (\ref{2.9}) gives
\begin{equation}
\Psi'_\gamma(\gamma t) \geq \mu(1) \Psi_\gamma(\gamma t)\psi(t)\geq
\mu(1) \Psi_\gamma(\gamma )\psi(t)\label{2.10}
\end{equation}
for all $t\in [1,+\infty[$. On the other hand from (\ref{2.0}) and
(\ref{2.1}) we
deduce
\begin{equation}
\|\varphi\|_{L^1}e^{\alpha \|\varphi\|_{L^1}}
{1\over b(t)}\geq \int_0^t b(s) \varphi(s)\, ds\geq 1 \label
{2.11}
\end{equation}
for all $t\in [1,+\infty[$. Finally since $\mu$ is increasing there exists
$\xi_0\geq 1$ such that
\begin{equation}
\mu({1\over \xi^2})\leq {1 \over 4\Lambda^2_0\|\psi\|_\infty
(C_0+\|\psi\|_\infty\mu(1))} Ê Ê\label{2.12}
\end{equation}
for all $\xi$ with $|\xi|\geq \xi_0$. Moreover Ê$\lim_{\gamma\to +\infty}
\Psi_\gamma(\gamma)=+\infty$ and then there exists $\gamma_0>0$ such that
\begin{equation}
\mu(1) \gamma \Psi_\gamma(\gamma) \int_0^1\varphi(s)\, ds
\geq (C_0 +\|\psi\|_\infty\mu({1\over \xi^2_0})) \mu({1\over\xi^2_0})\,\xi_0^4
\label {2.13}
\end{equation}
for all $\gamma\geq \gamma_0$. It is not restrictive to suppose also that
\begin{equation} \xi_0\geq  2\Lambda_0 \|\varphi\|_{L^1}e^{\alpha
\|\varphi\|_{L^1}}\quad {\rm and}\quad
\gamma_0\geq 4\Lambda^2_0 \|\varphi\|^2_{L^1}e^{2\alpha
\|\varphi\|_{L^1}}(C_0+\|\psi\|_\infty \mu(1)).
\label {2.14}
\end{equation}
We set
$$
\varepsilon_\xi=
\left\{
\begin{array}{l}
\displaystyle{{1\over \xi_0^2}\quad {\rm if}\quad |\xi|\leq \xi_0,}\\[0.5 cm]
\displaystyle{{1\over \xi^2}\quad {\rm if}\quad |\xi|\geq \xi_0.}
\end{array}
\right.
$$

Suppose first $|\xi|\leq \xi_0$. From (\ref{2.10}), (\ref{2.11}) and (\ref{2.13}) we
have
$$
\begin{array}{l}
\displaystyle{\gamma \Psi'_\gamma(\gamma t)\int_0^t b(s)
\varphi(s)\, ds}\\[0.3 cm]
\displaystyle{\qquad \geq \gamma Ê\mu(1) \Psi_\gamma(\gamma)\psi(t) b(t)
\int_0^t
\varphi(s)\, ds }\\[0.3 cm]
\displaystyle{\qquad\qquad \geq Êb(t) \psi(t) (C_0+\|\psi\|_\infty\mu({1\over
\xi_0^2}))
\mu({1\over \xi_0^2}) \xi_0^4}\\[0.3 cm]
\displaystyle{\qquad\qquad\qquad \geq Êb(t) \psi(t)(C_0
{\mu(\varepsilon_\xi)\over
\varepsilon_\xi}\xi^2 +\psi(t) \mu^2(\varepsilon_\xi)\xi^4)}
\end{array}
$$
for all $\gamma\geq \gamma_0$ and for all
$t\in [1,+\infty[$. Consequently
$$
\begin{array}{ll}
\displaystyle{\int_1^{+\infty}\int_{|\xi|\leq \xi_0}
(b(t)(a(t)\xi^2-\Phi'_\gamma(t))^2 +\gamma\Psi_\gamma'(\gamma t)\int_0^t
b(s)\varphi(s)\,ds) |\hat z(t,\xi)|^2\, d\xi dt}\\[0.5 cm]
\qquad\ \displaystyle{-\int_1^{+\infty}\int_{|\xi|\leq \xi_0} b(t)
\psi(t)(C_0 {\mu(\varepsilon_\xi)\over \varepsilon_\xi}\xi^2+ \psi(t)
\mu^2(\varepsilon_\xi)\xi^4) |\hat z(t,\xi)|^2\,d\xi dt \geq 0}
\end{array}
$$
for all $\gamma\geq\gamma_0$.

Suppose now $|\xi|\geq \xi_0$. If $a(t)\xi^2\geq 2\Phi'_\gamma(t)$ then
$$
b(t)(a(t)\xi^2-\Phi'_\gamma(t))^2\geq b(t){a^2(t)\over 4}\xi^4\geq
b(t){1\over 4\Lambda_0^2}\xi^4.
$$
As a consequence, from (\ref{2.12}), we have that
\begin{equation}
\begin{array}{l}
\displaystyle{b(t)\psi(t)(C_0{\mu(\varepsilon_\xi)\over \varepsilon_\xi}\xi^2
+\psi(t)\mu^2(\varepsilon_\xi) \xi^4)}\\[0.3 cm]
\displaystyle{\qquad = b(t)\psi(t)(C_0\mu({1\over \xi^2})
\xi^4 +\psi(t)\mu^2({1\over \xi^2}) \xi^4)}\\[0.3 cm]
\displaystyle{\qquad\qquad \leq
b(t)\|\psi\|_\infty(C_0+\|\psi\|_\infty\mu(1))\mu({1\over \xi^2})
\xi^4}\\[0.3 cm]
\displaystyle{\qquad\qquad\qquad\leq b(t){1\over 4\Lambda_0^2}\xi^4\leq
b(t)(a(t)\xi^2-\Phi'_\gamma(t))^2.}\\
\end{array}\label{2.15}
\end{equation}
If $a(t)\xi^2\leq 2\Phi'_\gamma(t)$ then (\ref{1}), (\ref{2.4})
and (2.12) imply that $$ \Psi_\gamma(\gamma t)\geq {\xi^2\over
2\Lambda_0\|\varphi\|_{L^1}e^{\alpha \|\varphi\|_{L^1}}}. $$ From
(\ref{2.9}) we infer $$
\begin{array}{l}
\displaystyle{\Psi'_\gamma(\gamma t)\geq {\xi^4\over
4\Lambda^2_0\|\varphi\|^2_{L^1}e^{2\alpha
\|\varphi\|_{L^1}}}\mu({2\Lambda_0\|\varphi\|_{L^1}e^{\alpha
\|\varphi\|_{L^1}}\over \xi^2})\psi(t)}\\[0.5 cm]
\qquad\qquad\quad\displaystyle{\geq
{\xi^4\over 4\Lambda^2_0\|\varphi\|^2_{L^1}e^{2\alpha
\|\varphi\|_{L^1}}}\mu(1/\xi^2)\psi(t).}
\end{array}
$$
Then
\begin{equation}
\gamma \Psi'_\gamma(\gamma t)\int_0^t b(s) \varphi(s)\, ds\geq
b(t) \psi(t) (C_0{\mu(\varepsilon_\xi)\over \varepsilon_\xi}\xi^2+\psi(t)
\mu^2(\varepsilon_\xi)\xi^4)
\label{2.16}
\end{equation}
for all $\gamma\geq \gamma_0$. Finally, (\ref{2.15}) and
(\ref{2.16}) give $$
\begin{array}{ll}
\displaystyle{\int_1^{+\infty}\int_{|\xi|\geq \xi_0}
(b(t)(a(t)\xi^2-\Phi'_\gamma(t))^2 +\gamma\Psi_\gamma'(\gamma t)\int_0^t
b(s)\varphi(s)\,ds) |\hat z(t,\xi)|^2\, d\xi dt}\\[0.5 cm]
\qquad\ \displaystyle{-\int_1^{+\infty}\int_{|\xi|\geq \xi_0} b(t)
\psi(t)(C_0 {\mu(\varepsilon_\xi)\over \varepsilon_\xi}\xi^2+ \psi(t)
\mu^2(\varepsilon_\xi)\xi^4) |\hat z(t,\xi)|^2\,d\xi dt \geq 0}
\end{array}
$$
for all $\gamma\geq\gamma_0$. The proof of Lemma~\ref{l1} is complete.


\begin{thebibliography}{10}

\bibitem{AN}
S. Agmon and L. Nirenberg,
Lower bounds and uniqueness theorems for solutions of differential equations
 in a Hilbert space,
Comm. Pure Appl. Math. {\bf 20} (1967), 207--229.
\bibitem{CoLees}
P.J. Cohen and M. Lees,
Asymptotic decay of solutions of differential inequalities,
Pacific J. Math. {\bf 11} (1961), 1235--1249.
\bibitem{DSP} D. Del Santo and M. Prizzi,  Backward uniqueness for parabolic operators whose
coefficients are non-Lipschitz continuous in time, J. Math. Pures Appl., to appear.
\bibitem{lax}
P.D. Lax, A stability theorem for solutions of abstract differential equations, and its
application to the study of the local behavior of solutions of elliptic equations,  Comm. Pure Appl.
Math.  {\bf 9}  (1956), 747--766.
\bibitem{lees}
M. Lees, Asymptotic behaviour of solutions of parabolic differential inequalities,  Canad. J.
Math.  {\bf 14}  (1962), 626--631.
\bibitem{LM} J. L. Lions and E. Magenes,
{\it Nonhomogeneous Boundary Value Problems and Applications I},
 Springer Verlag, Berlin Heidelberg New York 1972.
\bibitem{mill}
K. Miller, Nonunique continuation for uniformly parabolic and
elliptic equations in self-adjoint divergence form with H\"older
continuous coefficients,  Arch. Rational Mech. Anal. {\bf 54}
(1974), 105--117.
\bibitem{oga}
H. Ogawa, Lower bounds for solutions of parabolic differential
inequalities,  Canad. J. Math. {\bf 19}  (1967), 667--672.
\bibitem{prot}
M.H. Protter, Properties of solutions of parabolic equations and inequalities,  Canad. J.
Math.  {\bf 13}  (1961), 331--345.
\bibitem{tara}
S. Tarama, Local uniqueness in the Cauchy problem for second order elliptic equations with
non-Lipschitzian coefficients,  Publ. Res. Inst. Math. Sci.  {\bf 33}  (1997),  no. 1, 167--188.






\end{thebibliography}
\end{document}